\documentclass[leqno,12pt]{article}
\usepackage{latexsym,amsmath,amsfonts,amssymb,mathabx}
\usepackage{enumitem}
\usepackage{hyperref}
\usepackage{cancel}
\hypersetup{
    colorlinks,
    linkcolor={blue},
    citecolor={red},
    urlcolor={blue}
}

\usepackage[usenames,dvipsnames]{color}

\scrollmode

\newtheorem{thm}{Theorem}[section]
\newtheorem{lem}[thm]{Lemma}

\newtheorem{rem}[thm]{Remark}
\newtheorem{rems}[thm]{Remarks}

\newtheorem{df}[thm]{Definition}
\newtheorem{dfs}[thm]{Definitions}
\newtheorem{ex}[thm]{Example}

\newtheorem{as}[thm]{Assumption}

\newcommand{\wc}{\widecheck}
\newcommand{\wh}{\widehat}
\newcommand{\wt}{\widetilde}

\newcommand{\vS}{\varSigma}
\newcommand{\vO}{\varOmega}
\newcommand{\vT}{\varTheta}
\newcommand{\ov}{\overline}
\newcommand{\mf}{\mathfrak}
\newcommand{\vY}{\varUpsilon}

\newcommand{\uph}{\upharpoonright}

\def\N{{\mathbb N}}
\def\R{{\mathbb R}}

\def\F{{\mathcal F}}
\def\B{{\mathfrak B}}

\newcommand{\E}{\mathbb{E}}

\newcommand{\cnp}{\{N_t\}_{t\in\mathbb{R}_{+}}}
\newcommand{\cnpw}{counting process }

\newcommand{\clap}{\{T_n\}_{n\in\N_0}}

\newcommand{\cip}{\{W_n\}_{n\in\N}}

\newcommand{\abc}[1]{\textcolor{Green}{#1}}

\title{On the equivalence of various definitions\\of mixed Poisson processes}
\author{D.P. Lyberopoulos
\thanks{The author is indebted to the Public Benefit Foundation \textsc{Alexander S. Onassis}, which supported this research, under the Programme of Scholarships for Hellenes.}\;, 
 N.D. Macheras and S.M. Tzaninis}

\date{\today}

\begin{document}
\maketitle

\begin{abstract}
Under mild assumptions the equivalence of the mixed Poisson process with mixing parameter a real-valued random variable  to the one  with mixing distribution  as well as to the mixed Poisson process in the sense of Huang is obtained, and a characterization of each one of the above mixed Poisson processes in terms of disintegrations is provided. Moreover, some examples of ``canonical" probability spaces admitting counting processes satisfying the equivalence of all above statements are given.

 Finally, it is shown that our  assumptions are essential for the characterization of mixed Poisson processes in terms of disintegrations. \medskip

 \smallskip
\smallskip 

\smallskip

\end{abstract}


\section*{Introduction}\label{intro}

To the best of our knowledge, given a probability space $(\vO,\vS,P)$, there are five definitions for mixed Poisson processes (MPPs for short). The first one, involving birth processes, traces back to  Lundberg (cf. e.g. \cite{gr}, page 61), while the second one is that of the \emph {standard MPP with parameter a positive real-valued random variable} (cf. e.g. \cite{gr}, Definition 4.3). Another definition of MPPs associated with a family $\{P_{\wt{y}}\}_{\wt{y}\in\wt\vY}$  of probability measures on $\vS$ and with a probability measure $\nu$  on the $\sigma$-algebra $\sigma(\{P_{\bullet}(E):E\in\vS\})$ (written MPP$(\{P_{\wt{y}}\}_{\wt{y}\in\wt\vY},\nu)$ for short)  is due to Huang \cite{hu}, and is given in terms of interarrival processes, see Definition \ref{mpp} (b). The other two definitions refer to the cases of a MPP with mixing parameter a real-valued random variable $\vT$ (MPP$(\vT)$ for short) and a MPP with mixing distribution $U$ 
(MPP$(U)$ for short), see Definitions \ref{mpp} (a) and (c), respectively.

The equivalence of Lundberg's definition to that of the MPP$(U)$ is due to P. Albrecht \cite{Al}, while the equivalence of the definition of the standard mixed Poisson process with parameter a positive real-valued random variable  to that of the MPP$(\vT)$ is due to  R. F. Serfozo (see \cite{se1}, page 290 together with \cite{se2}, Theorem 3.1). 

In this paper we first investigate the equivalence of the definitions of a MPP$(\vT)$, a MPP$(U)$ and a MPP$(\{P_{\wt{y}}\}_{\wt{y}\in\wt\vY},\nu)$. It is easy to see that a MPP$(\vT)$ 
is always a MPP$(U)$. But the inverse implication does not seem to be in general true, as it is not always possible (given a MPP$(U)$) to construct a real-valued random variable $\vT$ such that $P_{\vT}=U$, and a disintegration of $P$ over $U$ consistent with $\vT$.
 
In Section \ref{mainr}  we prove  that given a counting process $N$, under  Assumption \ref{as1} and under the assumption that $P$ is perfect and $\vS$ is countably generated, the existence of a real-valued random variable $\wc{\vT}$ such that $N$ is a MPP$(\wc\vT)$ is equivalent with the existence of a distribution $U$ on $\B$ such that $N$ is a MPP$(U)$, and the latter is equivalent with the existence of a family $\{P_{\wt{y}}\}_{\wt{y}\in\wt\vY}$  of probability measures on $\vS$ and of a probability measure $\nu$  on $\sigma(\{P_{\bullet}(E):E\in\vS\})$ such that $N$ is a MPP$(\{P_{\wt{y}}\}_{\wt{y}\in\wt\vY},\nu)$, see Theorem \ref{thm1}. In Theorem \ref{thm2} we prove that under  Assumption \ref{as} and the existence of a disintegration of $P$ over $P_\vT$ consistent with a given real-valued random variable $\vT$, it follows that $N$ is a MPP$(\wh\vT)$ if and only if it is a MPP$(P_{\wh\vT})$ if and only if it is a MPP$(\{Q_{\wh{\theta}}\}_{\wh{\theta}\in\R},P_{\wh\vT})$, where $\wh\vT$ is a proper measurable function of $\vT$ and $\{Q_{\wh{\theta}}\}_{\wh{\theta}\in\R}$ is a proper family of probability measures on $\vS$. 

The proofs of Theorems \ref{thm1} and \ref{thm2} rely on two earlier results. The first one is due to  Lyberopoulos and Macheras where it is proven that under the existence of an appropriate disintegration of $P$ over $P_\vT$ a MPP$(\vT)$ can be reduced to an ordinary Poisson process under the disintegrating measures (see \cite{lm1v}, Proposition 4.4). The second one is due to Macheras and Tzaninis where it is proven  that under Assumption \ref{as} within the class of mixed renewal processes, a counting process is a MPP$(\vT)$ if and only if it has the $P$-Markov property (see \cite{lmt1n}, Theorem 2.11). For the definition of the $P$-Markov property we refer to e.g. \cite{Sc}, page 44.

 In Section \ref{exs} we provide two examples of ``canonical'' probability spaces where all assumptions of Theorems \ref{thm1} and  \ref{thm2} are valid. In particular, in both examples each of the assertions of Theorem \ref{thm2} is valid.

Finally, in Section \ref{cexs} we construct two counter-examples of non-trivial probability spaces where the characterization of MPPs in terms of disintegrations fails.


\section{Preliminaries}\label{prel}

By $\N$ is denoted the set of all natural numbers and $\N_0:=\N\cup\{0\}$.  The symbol $\R$ stands for the set of all real numbers, while $\overline{\R}:=\R\cup\{-\infty,+\infty\}$ and $\R^d$ denotes the Euclidean space of dimension $d\in\N$. Given a subset $A$ of a set $\vO$  we denote by $A^c$ the complement $\vO\setminus A$ of $A$ and by $\chi_A$ the indicator function of $A$. For a map $f:D\longmapsto E$ we denote by $R_f$ or by $f(D)$  the set $\{f(x):x\;\in D\}$, and for a set $A\subseteq D$ we denote by $f\upharpoonright A$ the restriction of $f$ to $A$, and by $f(A)$ the set $\{f(x):x\;\in A\}$.

Given a probability space $(\vO,\vS,P)$ a set $N\in\vS$ with $P(N)=0$ is called a $P$-{\bf null set}. For any two sets $A,B\in\vS$ we write $A=_P B$ if $P(A\triangle B)=0$. Given a measurable space $(\vY,H)$, for any two $\vS$-$H$-measurable maps $X,Y:\vO\longmapsto\vY$ we write $X=Y$ $P$-a.s. if $\{X\neq Y\}$ is a $P$-null set.

Given a topology  $\mf{T}$ on $\vO$ write ${\mf B}(\vO)$ for its {\bf Borel $\sigma$-algebra} on $\vO$, i.e. the $\sigma$-algebra generated by $\mf{T}$ and ${\mf B}:={\mf B}(\R)$, ${\overline{\mf B}}:={\mf B}(\overline{\R})$, ${\mf B}_d:={\mf B}(\R^d)$  and $\mf{B}_{\N}:={\mf B}(\R^{\N})$ for the Borel $\sigma$-algebra of subsets of $\R$, $\overline{\R}$, $\R^d$ and $\R^{\N}$ under the corresponding Euclidean topologies, respectively, while ${\mathcal{L}}^{1}(P)$ stands for the family of all real-valued $P$-integrable functions on $\vO$. Functions that are $P$-a.s. equal are not identified.

For the definitions of {\bf real-valued random variables}, {\bf random variables} and {\bf random vectors} we refer to Cohn \cite{co}, pages 308 and 318.

Given two probability spaces $(\vO,\vS,P)$ and $(\vY,H,Q)$ as well as a $\vS$-$H$-measurable map $X:\vO\longmapsto\vY$ we denote by $\sigma(X):=\{X^{-1}(B): B\in H\}$ the $\sigma$-algebra generated by $X$, while $\sigma(\{X_i\}_{i\in I}):=\sigma\bigl(\bigcup_{i\in{I}}\sigma(X_i)\bigr)$ stands for the $\sigma$-algebra generated by a family $\{X_i\}_{i\in I}$ of $\vS$-$H$-measurable maps from $\vO$ into $\vY$.

For any $d$-dimensional  random vector $X$ on $\vO$ we apply the notation $P_{X}={\bf{K}}(\theta)$ in the meaning that $X$ is distributed according to the law ${\bf{K}}(\theta)$, where $\theta\in\R^d$. In particular, $\mathbf{P}(\theta)$ and $\mathbf{Exp}(\theta)$, where $\theta$ is positive parameter, stand for the law of Poisson and exponential distribution, respectively (cf. e.g. \cite{Sc}).

We write $\E[X|\F]$ for a conditional expectation of $X$ given $\F$ (see \cite{co}, page 342 for the definition). For $X:=\chi_E\in\mathcal{L}^1(P)$ with $E\in\vS$ we set $P(E\mid\mathcal{F}):=\E_P[\chi_E\mid\mathcal{F}]$.

Given a real-valued random variable $X$ on $\vO$ and a random vector $\vT:\vO\longmapsto\R^d$, a conditional distribution of $X$ over $\vT$ is a map $P_{X|\vT}$ from $\B\times\vO$ into $[0,1]$ such that 
\begin{description}
\item[(cd1)] for each $\omega\in\vO$ the set-function $P_{X|\vT}(\bullet,\omega)$ is a probability measure on $\B$;
\item[(cd2)] for each $B\in\B$ we have
$$
P_{X|\vT}(B,\bullet)=P(\vT^{-1}(B)\mid\sigma(\vT))\quad P\upharpoonright\sigma(\vT)\text{-a.s.,}
$$
where $P_{X|\vT}(B,\bullet)$ is $\sigma(\vT)$-measurable for any fixed $B\in\B$.
\end{description}
For simplicity we write $k:=P_{X|\vT}$ and define the map $K(\vT)$ from $\B\times\vO$ into $[0,1]$ by means of 
$$
K(\vT)(B,\omega):=(k(B,\bullet)\circ\vT)(\omega)\quad \forall\,B\in\B \,\,\,\,\forall\,\omega\in\vO.
$$
Then for $\theta=\vT(\omega)$ with $\omega\in\vO$ the probability measures $k(\bullet,\theta)$ are distributions on $\B$ and so we may write ${\bf K}(\theta)(\bullet)$ instead of $k(\bullet,\theta)$. Consequently, in this case $K(\vT)$ will be written by ${\bf{K}}(\vT)$. 

For any real-valued random variables $X$, $Y$ on $\vO$ we say that $P_{X|\vT}$ and $P_{Y|\vT}$ are $P\upharpoonright\sigma(\vT)$-equivalent and we write $P_{X|\vT}=P_{Y|\vT}$ $P\upharpoonright\sigma(\vT)$-a.s., if there exists a $P$-null set $N\in\sigma(\vT)$ such that for any $\omega\notin N$ and $B\in \B$ the equality $P_{X|\vT}(B,\omega)=P_{Y|\vT}(B,\omega)$ holds true.\smallskip

{\em{From now on $(\vO,\vS,P)$ is a fixed probability space, while $(\vY,H):=(\R,\B)$, $(\varXi,Z):=(\R^d,\B_d)$.}}


\section{The results}\label{mainr}
We first recall some additional background material, needed in this section.

A family $N:=\cnp$ of  random variables from $(\vO,\vS)$ into $(\ov{\R},\ov{\B})$ is called a {\bf counting process} if there exists a $P$-null set $\vO_N\in\vS$ such that the process $N$ restricted on $\vO\setminus\vO_N$ takes values in $\N_0\cup\{\infty\}$, has right-continuous paths, presents jumps of size (at most) one, vanishes at $t=0$ and increases to infinity. Denote by $T:=\clap$ and $W:=\cip$  the {\bf{arrival process}} and {\bf{interarrival process}} respectively (cf. e.g. \cite{Sc}, Section 1.1, page 6 for the definition) associated with $N$. \smallskip

Recall that for a random vector $\vT:\vO\longmapsto \R^d$ a family $\{X_i\}_{i\in I}$ of real-valued random variables  $X_i$ on  $\vO$ 
\begin{itemize}
\item[$\bullet$] is {\bf $P$-conditionally (stochastically) independent given} $\vT$, if for each $n\in\N$ with $n\geq 2$ we have
$$
P(\bigcap^{n}_{j=1} \{X_{i_j}\leq x_{i_j}\}\mid\sigma(\vT))=\prod^{n}_{k=1} P(\{X_{i_j}\leq x_{i_j}\}\mid\sigma(\vT))\qquad\quad P\upharpoonright\sigma(\vT)-\mbox{a.s.}
$$
whenever $i_1,\ldots,i_n$ are distinct members of $I$ and $(x_{i_1},\ldots,x_{i_n})\in\R^n$;
\item[$\bullet$] is {\bf $P$-conditionally identically distributed given} $\vT$, if 
$$
P\bigl(F\cap X_i^{-1}(B)\bigr)=P\bigl(F\cap X_j^{-1}(B)\bigr)
$$
whenever $i,j\in I$, $F\in\sigma(\vT)$ and $B\in \B$.
\end{itemize}

{\em For the rest of the paper we simply write ``conditionally'' in the place of ``conditionally given $\vT$'' whenever $\vT$ is clear from the context.}

\begin{df}\label{rcp}
\normalfont
Let $Q$ be a probability measure on $\B_d$. A family $\{P_\theta\}_{\theta\in\R^d}$ of probability measures on $\vS$ 
is called a {\bf disintegration} of $P$ over $Q$ if
\begin{description}
\item[(d1)] 
for each $D\in\vS$ the map $\theta\longmapsto P_\theta(D)$ is $\B_d$-measurable;
\item[(d2)]
$\int P_{\theta}(D)\,Q(d\theta)=P(D)$ for each $D\in\vS$.
\end{description}
If $\vT:\vO\longmapsto\R^d$ 
is an inverse-measure-preserving random vector (i.e. $P_\vT(B)=Q(B)$ for each $B\in \B_d$), a disintegration $\{P_{\theta}\}_{\theta\in\R^d}$ of $P$ over $Q$ is called {\bf consistent} with $\vT$ if, for each $B\in \B_d$, the equality $P_{\theta}(\vT^{-1}(B))=1$ holds for $Q$-almost every $\theta\in B$.
\end{df}

{\em Remark.} 
If $\vS$ is countably generated (cf. e.g. \cite{co}, Section 3.4, page 102 for the definition) and $P$ is perfect (see \cite{fa}, p. 291 for the definition), then there always exists a disintegration $\{P_{\theta}\}_{\theta\in\R^d}$ of $P$ over $Q$ consistent with any inverse-measure-preserving 
random vector  $\vT:\vO\longmapsto\R^d$ (see \cite{fa}, Theorems 6 and 3) and this means, that in most cases appearing in applications (e.g. Polish spaces) disintegrations
consistent with a random vector exist.
\smallskip

{\em Throughout what follows, unless stated otherwise, $N:=\cnp$ is a counting process, $T:=\left\{T_{n}\right\}_{n\in\N_0}$ is an  arrival process, $W:=\left\{W_{n}\right\}_{n\in\N}$ is its induced  interarrival process  and without loss of generality we may and do assume that $\vO_N=\emptyset$.}\smallskip

A Poisson process $N$ with respect to $P$ with parameter $\theta>0$ is denoted by $P$-PP$(\theta)$.

\begin{dfs}\label{mpp}
\normalfont
A \cnpw $N$ is:

{\bf{(a)}} a {\bf{mixed Poisson process with  mixing parameter a real-valued random variable $\vT$}} such that $P_{\vT}\bigl(\left(0,\infty\right)\bigr)=1$ (written $P$-MPP$(\vT)$ for short), if it has $P$-conditionally  independent and $P$-conditionally  stationary increments (cf. e.g. \cite{Sc}, Section 4.1, page 86 for the definition)
and condition 
$$ 
\forall\; t\in (0,\infty)\qquad [P_{N_{t}\mid\vT}={\bf{P}}\left(\vT t\right)\quad P\upharpoonright\sigma(\vT)-\text{a.s.}]
$$
holds true (cf. e.g. \cite{Sc}, page 87);

{\bf{(b)}} a {\bf mixed Poisson process associated with $\{P_{\wt{y}}\}_{\wt{y}\in\wt\vY}$ and $\nu$} ($P$-MPP$(\{P_{\wt{y}}\}_{\wt{y}\in\wt\vY},\nu)$ for short), if for every $r\in\N$ and for 
all $w_1,\ldots,w_r>0$ condition
$$ P\Bigl(\bigcap_{k=1}^{r}\{W_k\leq w_k\}\Bigr) =\int\prod_{k=1}^{r}P_{\wt{y}}(\{W_k\leq w_k\})\,\nu(d\wt{y}),$$
holds true, where $\{P_{\wt{y}}\}_{\wt{y}\in\wt\vY}$ is a family of probability measures on $\vS$ and $\nu$ is a probability measure on $B(\vS):=B(\wt\vY,\vS):=\sigma(\{P_{\bullet}(E):E\in\vS\})$ such that $W$ is $P_{\wt y}$-independent and $\left(P_{\wt{y}}\right)_{W_{n}}={\bf{Exp}}\left(\alpha(\wt{y})\right)$  
 for every $n\in\mathbb N$ and for $\nu$-a.a. $\wt{y}\in\wt{\varUpsilon}$, where
$\alpha$ is a positive measurable function on $\R$ (see \cite{hu}, page 2);

{\bf{(c)}} a {\bf{mixed Poisson process with mixing distribution}} $U:\mathfrak B\left((0,\infty)\right)\longmapsto [0,1]$ (written $P$-MPP$(U)$ for short) if
$$
P\left(\bigcap^{m}_{j=1}\{N_{t_{j}}-N_{t_{j-1}}=\kappa_{j}\}\right)
=\int_{(0,\infty)}\prod^{m}_{j=1} e^{-\theta(t_{j}-t_{j-1})}\frac{\big(\theta(t_{j}-t_{j-1})\big)^{\kappa_{j}}}{\kappa_{j}!}\,U(d\theta)
$$
holds for all $m \in\mathbb N$ and $t_{0}, t_{1},\ldots,t_{m}\in\mathbb R_{+}$  with $0 =t_{0}< t_{1}<\ldots<t_{m}$ and for all 
$\kappa_{j}\in\mathbb N_{0}$, $j\in\{1,\ldots, m\}$ (cf. e.g. \cite{SZ}, page 9).
\end{dfs}

The following definition has been introduced in \cite{lmt1n}, Definitions 2.3.

\begin{df}\label{rd}
\normalfont
A counting process $N$ is called an  {\bf extended MRP  with mixing parameters $\vT$ and $h$, and interarrival time conditional distribution ${\bf K}(h(\vT))$} (written $P$-eMRP$({\bf K}(h(\vT)))$ for short), if $h$ is a $\R^k$-valued $\mathfrak B(D)$-$\mathfrak B_k$-measurable function on $D\in\B_d$ with $R_\vT\subseteq D$ for $k\in\N$, if the induced  interarrival process $W$ is $P$-conditionally independent and 
$$
\forall\;n\in\N\qquad [P_{W_n|\vT}={\bf{K}}\left(h(\vT)\right)\qquad P\upharpoonright\sigma(\vT)-\text{a.s.}].
$$
In particular, if $k=d$ and $h=id_{\R^d}$ then $N$ is  a {\bf $P$-MRP with  interarrival time distribution $\mathbf{K}(\vT)$} (written (written $P$-MRP$({\bf{K}}(\vT))$ for short). Moreover if $h=id_{\R^d}$ and if there exists a $\theta_0\in\R^d$ with $P(\{\vT=\theta_0\})=1$ then $N$ is a {\bf $P$-renewal process with  interarrival time distribution $\mathbf{K}(\theta_0)$} (written $P$-RP$(\mathbf{K}(\theta_0))$ for short).

Without loss of generality we may and do assume that 
\begin{equation}\label{rd1}
\forall\, n\in\N\qquad [P_{W_n|\vT}={\bf{K}}\left(h(\vT)\right)].
\end{equation}
\end{df}

{\em From now on, unless stated otherwise, $\vT$ is a positive real-valued random variable on $\vO$.}\smallskip

The following assumption is a special case of Assumption 2.6 from \cite{lmt1n}.

\begin{as}\label{as}
\normalfont
Let $D\in\B$ with $R_\vT\subseteq D$, $h:D\longmapsto\R$ be a $\B(D)$-measurable function, let $N$ be a $P$-eMRP$({\bf K}(h(\vT)))$ and let $\{P_\theta\}_{\theta\in D}$  be  a disintegration of $P$ over $P_{\vT}$ consistent with $\vT$. It follows by \cite{lm6z3}, Lemma 3.5 together with condition (\ref{rd1}) that 
\begin{equation}\label{rd2}
\forall\, n\in\N\quad\forall\,\theta\in{D}\qquad [(P_{\theta})_{W_n}={\bf{K}}\left(h(\theta)\right)].
\end{equation}
For any $\theta\in D$ and $t\in\R_+$ put
$$
F_{h(\theta)}(t):=P_\theta(\{W_n\leq t\})\quad\text{for all}\quad n\in\N.
$$
Clearly the function $F_{h(\theta)}$ depends on the distribution of $W_n$ and, because of condition (\ref{rd2}), on $h$. We say that $N$, $h$ and $\{P_\theta\}_{\theta\in D}$  satisfy Assumption \ref{as}, if there exists a $P_\vT$-null set $L_h:=L_{h,N,\{P_\theta\}_{\theta\in D}}$ in $\mathfrak B(D)$ such that 
for any $\theta\notin L_{h}$ the function $F_{h(\theta)}$ is continuously differentiable on $(0,\infty)$, there exists a function $C\in\mathcal L^1(P_{h(\vT)})$ with $0<F'_{h(\theta)}(t)<C(h(\theta))$ for each $t>0$, and the function $p_h:D\setminus L_h\longmapsto\R$ defined by means of $p_h(\theta):=p_{h,1}(\theta):=\lim_{t\rightarrow 0} F'_{h(\theta)}(t)$ is positive and injective.

For the special case  $D=\R$ and $h:=id_{\R}$ we write for simplicity $L$, $F_\theta$ and $p_1$ in the place of $L_h$, $F_{h(\theta)}$ and $p_{h}$ respectively, and we say that $N$ and $\{P_\theta\}_{\theta\in\R}$  satisfy Assumption \ref{as}.
\end{as}

\begin{as}\label{as1}
\normalfont
Given a counting process $N$ there exists a real-valued random variable $\vT$ on $\vO$ and disintegration $\{P_{\theta}\}_{\theta\in D}$ of $P$ over $P_\vT$ consistent with $\vT$ satisfying together with $N$ Assumption \ref{as}.
\end{as}

\begin{thm}\label{thm1}
For a counting process $N$ consider the following assertions:
\begin{enumerate}
\item
there exists a  real-valued random variable $\wc{\vT}$ on $\vO$ such that $N$ is a $P$-MPP$ (\wc{\vT})$ with respect to $P$;
\item
there exists a family $\{P_{\wt{y}}\}_{\wt{y}\in\wt{\varUpsilon}}$ of probability measures on $\vS$ and a probability measure $\nu$ on $B(\vS)$ such that $N$ is a $P$-MPP$(\{P_{\wt{y}}\}_{\wt{y}\in\wt\vY},\nu)$;
\item
there exists a  real-valued random variable $\wc{\vT}$ on $\vO$  and a disintegration $\{Q_{\wc{\theta}}\}_{\wc{\theta}\in\R}$ of $P$ over $P_{\wc{\vT}}$ consistent with $\wc{\vT}$ such that $N$ is a 
$Q_{\wc{\theta}}$-PP$(\wc{\theta})$  for $P_{\wc{\vT}}$-a.a. $\wc{\theta}\in\R$;
\item
there exists a distribution $U:\B\longmapsto[0,1]$ with $U((0,\infty))=1$ such that $N$ is a $P$-MPP$(U)$.
\end{enumerate}
Then $(iii)\Longrightarrow{(i)}\Longrightarrow{(iv)}$.

Moreover, if $P$ is perfect and $\vS$ is countably generated then statements $(i)$ and $(iii)$ are equivalent
and $(i)\Longrightarrow(ii)$. 

If in addition, Assumption \ref{as1} holds true, then all statements $(i)$ to $(iv)$ are equivalent.
\end{thm}

{\bf{Proof.}} 
First note that implication $(iii)\Longrightarrow(i)$ is immediate by Proposition 4.4 of \cite{lm1v}, while the implication $(i)\Longrightarrow(iv)$ follows by an easy computation.

Assume now that $P$ is perfect and $\vS$ is countably generated. Then the equivalence $(i)\Longleftrightarrow(iii)$ follows by \cite{lm1v}, Proposition 4.4, since under this assumption for any real-valued random variable $\wc{\vT}$ on $\vO$ there always exists a disintegration $\{Q_{\wc{\theta}}\}_{\wc{\theta}\in\R}$ of $P$ over $P_{\wc{\vT}}$ consistent with $\wc{\vT}$ (see \cite{fa}, Theorems 6 and 3). 

Ad $(i) \Longrightarrow (ii)$: If $(i)$ is true, since $(i)$ is equivalent with $(iii)$, it follows that there exists a disintegration $\{Q_{\wc{\theta}}\}_{\wc{\theta}\in\R}$ of $P$ over $P_{\wc{\vT}}$ consistent with $\wc{\vT}$ such that $N$ is a $Q_{\wc{\theta}}$-PP$(\wc{\theta})$  for $P_{\wc{\vT}}$-a.a. $\wc{\theta}\in\R$. Then applying \cite{lmt1n}, Proposition 2.2, we obtain $(ii)$ for $\{P_{\wt y}\}_{\wt y\in\wt\vY}:=\{Q_{\wc{\theta}}\}_{\wc{\theta}\in\R}$ and $\nu:=P_{\wc\vT}$.

Assume now in addition that Assumption \ref{as1} holds true.

Ad $(ii) \Longrightarrow (i)$: If assertion $(ii)$ holds true, then by Definition \ref{mpp} (b) together with e.g. \cite{Sc}, Theorem  2.3.4, we obtain that $N$ is a $P_{\wt y}$-PP$(\alpha(\wt{y}))$ for $\nu$-a.a. $\wt{y}\in\wt{\varUpsilon}$; hence applying \cite{Sc}, Lemma 2.3.1 we deduce that for $\nu$-a.a. $\wt y\in\wt \vY$ the equality 
\begin{equation}
\label{1}
P_{\wt y}(\bigcap_{j=1}^{m}\{N_{t_j}-N_{t_{j-1}}=n_j\})=\frac{n!}{\prod_{j=1}^m n_j!}\cdot\prod_{j=1}^m\Big(\frac{t_j-t_{j-1}}{t_m}\Big)^{n_j}\cdot P_{\wt y}(\{N_{t_m}=n\})
\end{equation}
holds true for any $m\in\N$, any $t_0,t_1,\ldots,t_m \in\R_+$ such that $0=t_0<t_1<\cdots<t_m$ and any $n_1,\ldots,n_m\in\N_0$ such that $\sum_{j=1}^{m} n_j=n$. The latter implies again for 
 $\nu$-a.a. $\wt y\in\wt \vY$ the equality 
\begin{equation}
\label{2}
P_{\wt y}(\{N_{s}=k\}\cap\{N_{t}-N_{s}=n-k\})={n\choose k}\cdot\Big(\frac{s}{t}\Big)^{k}\cdot\Big(1-\frac{s}{t}\Big)^{n-k}\cdot P_{\wt y}(\{N_t=n\})
\end{equation}
for all $s,t\in(0,\infty)$ such that $s<t$ and all $k,n\in\N_0$ such that $k\leq n$. Putting $\F_N:=\sigma(\{N_t\}_{t\in\R_+}$, $\F_W:=\sigma(\{W_n\}_{n\in\N}$ and $\F_T:=\sigma(\{T_n\}_{n\in\N_0}$ we  then get $\F_T=\F_W$ and $\F_T=\F_N$ (cf. e.g. \cite{Sc}, Lemmas 1.1.1 and 2.1.3 respectively).\medskip

\begin{it}
Claim. The family $\{P_{\wt y}\}_{\wt y\in\wt\vY}$ of probability measures is a disintegration of $P\uph\F_N$ over $\nu$.\smallskip

Proof.
\end{it} 
Since $\F_N=\F_W$, it is sufficient to show that $\{P_{\wt y}\}_{\wt y\in\wt\vY}$  is a disintegration of $P\uph\F_W$ over $\nu$. 

It follows by  Definition \ref{mpp} (b) that the property (d1) holds true; hence it is enough to show (d2) for any $E\in\F_W$. Put $\mathcal G:=\bigcup_{n\in\N}\sigma(W_n)$. Due to Definition \ref{mpp} (b) we have that (d2) is satisfied for each $\{W_n\leq w_n\}$ where $w_n>0$ for any $n\in\N$.

Denote by $\mathcal G_{\bigcap}$ be the generator of $\F_W$ consisting of $\mathcal G$ and all finite intersections of elements of $\mathcal G$ and put
$$
\mathcal D:=\{E\in\F_W : P(E)=\int P_{\wt y}(E)\,\nu(d\wt y)\}.
$$
Then it easy to prove that the family $\mathcal D$ is a Dynkin class  containing  $\mathcal G_{\bigcap}$; hence by the Monotone Class Theorem we get that $\mathcal D=\F_W$ and the above claim  follows.  \hfill$\Box$\smallskip

Fix on arbitrary $m\in\N$, $t_0,t_1,\ldots,t_{m+1} \in\R_+$ such that $0=t_0<t_1<\cdots<t_{m+1}$ and  $n_0,n_1,\ldots,n_{m+1}\in\N_0$ such that $0=n_0\leq n_1\leq\cdots\leq n_{m+1}$. Using  the above claim  and the equalities \eqref{1} and \eqref{2} we get by standard computations that
\begin{eqnarray*}
\lefteqn{P(\bigcap_{j=1}^{m}\{N_{t_j}=n_j\})\cdot P(\{N_{t_m}=n_m\}\cap\{N_{t_{m+1}}=n_{m+1}\})}&&\\
&\qquad\qquad\qquad\qquad =& P(\bigcap_{j=1}^{m+1}\{N_{t_j}=n_j\})\cdot P(\{N_{t_m}=n_m\}),
\end{eqnarray*}
or equivalently that $N$ has the $P$-Markov property. Since  $N$ has the $P$-Markov property and Assumption  \ref{as1} is valid, we may apply \cite{lmt1n}, Proposition 2.7, to obtain assertion $(i)$.

Ad $(iv) \Longrightarrow (i)$: If $(iv)$ is valid then by Theorem 4.2 from \cite{SZ}, $N$ has the $P$-Markov property. By Assumption \ref{as1} there exists  a real-valued random variable $\vT$ on $\vO$ and a disintegration $\{P_\theta\}_{\theta\in D}$ of $P$ over $P_\vT$ consistent with $\vT$ satisfying together with $N$ Assumption \ref{as}. Thus, we may apply \cite{lmt1n} Proposition 2.7 to get $(i)$. 
\hfill$\Box$

\begin{thm}\label{thm2}
Let $N$ be a $P$-eMRP$({\bf K}(h(\vT)))$ and let $\{P_{\theta}\}_{\theta\in D}$ be a disintegration of $P$ over $P_{\vT}$ consistent with $\vT$ satisfying together with $N$ and $h$ Assumption \ref{as}. Suppose  that there exists a $P_\vT$-null set $L_0\in\B(D)$ such that $h\uph D\setminus L_0$ is injective. Put $O_h:=L_0\cup L_h$ and $\wh{\vT}(\omega):=(p_h\circ\vT)(\omega)$ if $\omega\in\vT^{-1}(D\setminus O_h)$, where $L_h$ and $p_h$ are as in Assumption \ref{as}, and denote again by $\wh\vT$ any measurable extension of $\wh\vT$ from $\vT^{-1}(D\setminus O_h)$ to $\vO$. For any fixed $A\in\vS$ put
$$
Q_{\wh{\theta}}\left(A\right):=\begin{cases}(P_{\bullet}(A)\circ  p^{-1}_h)(\wh{\theta}) &\text{if}\;\; \wh{\theta}\in p_h(D\setminus O_h); \\ P\left(A\right)& \text{otherwise}.\quad\quad\end{cases}
$$ 
Then $\{Q_{\wh\theta}\}_{\wh\theta\in\R}$ is a disintegration of $P$ over $P_{\wh{\vT}}$ consistent with $\wh{\vT}$, and the following are equivalent:
 
\begin{enumerate}
\item
$N$ is a $P$-MPP$(\wh{\vT})$;
\item
$N$ is a $P$-MPP$(\{Q_{\wh\theta}\}_{\wh\theta\in\R},P_{\wh\vT})$;
\item
$N$ is a $Q_{\wh{\theta}}$-PP$(\wh{\theta})$  for $P_{\wh\vT}$-a.a. $\wh{\theta}\in \R$;
\item
$N$ is a $P$-MPP$(P_{\wh{\vT}})$.
\end{enumerate}
\end{thm}

{\bf{Proof.}} The fact that $\{Q_{\wh\theta}\}_{\wh\theta\in\R}$ is a disintegration of $P$ over $P_{\wh{\vT}}$ consistent with $\wh{\vT}$ is a consequence of \cite{lmt1n}, Lemma 2.4. 

The equivalence   $(i)\Longleftrightarrow(iii)$ is due to \cite{lm1v}, Proposition 4.4.

Ad $(i) \Longrightarrow (ii):$ Since assertion $(i)$ holds true and  $\{Q_{\wh{\theta}}\}_{\wh{\theta}\in\R}$ is a disintegration of $P$ over $P_{\wh{\vT}}$ consistent with $\wh{\vT}$, it follows by \cite{lmt1n}, Proposition 2.2 that $(ii)$ is valid.

Ad $(ii) \Longrightarrow (i):$ If $(ii)$ holds true, we get as in the proof of Theorem \ref{thm1}, $(ii) \Longrightarrow (i)$,  that $N$ has the $P$-Markov property; hence by \cite{lmt1n}, Theorem 2.11,  we obtain $(i)$.

The implication \textit{(i)}$\Longrightarrow$\textit{(iv)} follows by an easy computation.

Ad $(iv) \Longrightarrow (i):$ If $(iv)$ is valid then  $N$ has the $P$-Markov property (see \cite{SZ}, Theorem 4.2). So, by \cite{lmt1n}, Theorem 2.11, assertion $(i)$ follows.
\hfill$\Box$

\begin{rems}\label{rem2}
\normalfont
{\bf (a)} If the assumptions of Theorem \ref{thm2} are satisfied and if one of its assertions $(i)$-$(iv)$ is valid, then $p_h$ and $h$ coincide outside the $P_\vT$-null set $\wt L_3\cup O_h$, and $\wh\vT$ and $\wt\vT$ coincide outside of the $P$-null set $\vT^{-1}(O_h)$.\smallskip

In fact, let $\{Q_{\wh\theta}\}_{\wh\theta\in\R}$ be as in Theorem \ref{thm2} and assume that $(iii)$ holds true. It then follows that $(Q_{\wh\theta})_{W_n}={\bf Exp}(\wh\theta)$ and that $W$ is $Q_{\wh\theta}$-independent for any $\wh{\theta}\in p_h(D\setminus O_h)$ (cf. e.g. \cite{Sc}, Theorem 2.3.4). Applying now \cite{lmt1n}, Lemma 2.4, for $p_h(O_h)$ and $p_h$ in the place of $L_0$ and $h$ respectively, we obtain that $(P_{\theta})_{W_n}={\bf Exp}(p_h(\theta))$ and that $W$ is $P_{\theta}$-independent for any $\theta\in D\setminus O_h$. Since $N$ is $P$-eMRP$({\bf K}(h(\vT)))$, due to \cite{lmt1n}, Lemma 2.9, we can find a $P_\vT$-null set $\wt L_3\in\B(D)$ such that $(P_\theta)_{W_n}={\bf K}(h(\theta))$ for any $\theta\in D\setminus \wt L_3$ and $n\in\N$.

As a consequence, we deduce that for any $\theta\in D\setminus (\wt L_3\cup O_h)$ and $\wh\theta=p_h(\theta)$ conditions
$$
{\bf Exp}(\wh\theta)=(Q_{\wh\theta})_{W_n}=(P_{\theta})_{W_n}={\bf K}(h(\theta))
$$
hold true, implying that $p_h(\theta)=h(\theta)$ for any $\theta\in D\setminus(\wt L_3\cup O_h)$.\smallskip

{\bf (b)} It is worth noticing that if all statements $(i)$ to $(iv)$ of Theorem \ref{thm2} are equivalent and if one of them is valid, then its assumptions are necessary. 
More precisely, let $h$ and $L_0$ be as in Theorem \ref{thm2} such that $\E[h(\vT)]<\infty$, let $N$ be a counting process and $\wt\vT:=h\circ\vT$. Because of (a) we may take $h$ and $\wt\vT$ in the place of $p_h$ and $\wh\vT$ respectively. 
 Assume that assertions $(i)$-$(iv)$ of Theorem \ref{thm2} are all equivalent and each of them is valid with $\wt\vT$, $h$ and $L_0$ in the place of $\wh\vT$, $p_h$ and $O_h$. Then $N$ is  a $P$-eMRP$({\bf K}(h(\vT)))$, and there exists a disintegration $\{P_{\theta}\}_{\theta\in D}$  of $P$ over $P_{\vT}$ consistent with $\vT$ satisfying together with $N$ and $h$ Assumption \ref{as}.\smallskip

In fact, since $(iii)$ is valid, there exists a disintegration $\{Q_{\wt{\theta}}\}_{\wt{\theta}\in\R}$ of $P$ over $P_{\wt{\vT}}$ consistent with $\wt{\vT}$ such that the counting process $N$ is a 
PP$(\wt{\theta})$ with respect to $Q_{\wt{\theta}}$ for 
any $\wt\theta\in h(D\setminus L_0)$. 
The latter is equivalent with the fact that $(Q_{\wt\theta})_{W_n}={\bf Exp}(\wt\theta)$ and that $W$ is $Q_{\wt\theta}$-independent for any $\wt\theta\in h(D\setminus L_0)$ (cf. e.g. \cite{Sc}, Theorem 2.3.4). For any $\theta\in D$ and $A\in\vS$ define 
$$
P_{\theta}\left(A\right):=\begin{cases}(Q_{\bullet}(A)\circ h)(\theta) &\text{if}\;\; \theta\in D\setminus L_0; \\ P\left(A\right)& \text{otherwise}.\quad\quad\end{cases}
$$
Applying \cite{lmt1n}, Lemma 2.4, we obtain that $\{P_{\theta}\}_{\theta\in D}$ of $P$ over $P_{\vT}$ consistent with $\vT$ and that $(P_{\theta})_{W_n}={\bf Exp}(h(\theta))$ and $W$ is $P_{\theta}$-independent for any $\theta\in D\setminus L_0$.  Applying now \cite{lmt1n}, Lemma 2.9, together with \cite{lm6z3}, Lemma 3.6, we obtain that $P_{W_n|\vT}={\bf Exp}(h(\vT))$ $P\upharpoonright\sigma(\vT)$-a.s. and that $W$ is $P$-conditionally independent; hence $N$ is a $P$-eMRP$({\bf K}(h(\vT)))$.

It remains to show that $\{P_\theta\}_{\theta\in D}$, $N$ and $h$ satisfy Assumptions \ref{as}.

In fact, for any $\theta\in D\setminus L_0$, $t\in\R_+$ and  $n\in\N$ put $F_{h(\theta)}(t):=P_\theta(\{W_n\leq t\}):=1-e^{-h(\theta) t}$. Clearly, $F_{h(\theta)}$ is continuously differentiable on $(0,\infty)$. Define the map $C\in \mathcal L_1(P_{h(\vT)})$ by $C(h(\theta)):=h(\theta)$ for any $\theta\in D\setminus L_0$, and for any fixed $\theta\in D\setminus L_0$ define the density $f_{h(\theta)}:=F'_{h(\theta)}$ by $f_{h(\theta)}(t):=h(\theta)\cdot e^{-h(\theta) t}$ for any $t>0$.  Clearly, for any fixed $\theta\in D\setminus L_0$, the density $f_{h(\theta)}$ is dominated by $C(h(\theta))$, and the function $\lim_{t\rightarrow 0} f_{h(\theta)}(t)=h(\theta)$ is positive and injective; hence $\{P_\theta\}_{\theta\in D}$, $N$ and $h$ satisfy Assumptions \ref{as}.  
\end{rems}


\section{Examples}\label{exs}
By $(\vO\times\vY,\vS\otimes{H},P\otimes{Q})$ is denoted the product probability space of $(\vO,\vS,P)$ and $(\vY,H,Q)$, and by $\pi_{\vO}$ 
and $\pi_{\vY}$ the canonical projections from $\vO\times\vY$ onto $\vO$ and $\vY$, respectively.
 
In this section we first provide an example of ``canonical" probability spaces admitting extended MRPs. Next we present, as special cases,  two examples of probability spaces satisfying all assumptions of Theorems \ref{thm1} and \ref{thm2}. 
In particular, in both examples each of the assertions of Theorems \ref{thm1} and \ref{thm2} is valid.
\smallskip

{\em Throughout what follows, we put $\varUpsilon:=(0,\infty)$, $H:=\mathfrak B(\varUpsilon)$, $\wt{\vO}:=\varUpsilon^{\mathbb N}$, $\vO:=\wt{\vO}\times G$ for $G\in\B$, $\wt{\vS}:=\mathfrak B(\wt{\vO})$ and $\vS:=\mathfrak B(\vO)=\mathfrak B(\wt{\vO})\otimes \mathfrak B(G)$ for simplicity.}\smallskip

The next example is a special case of Example 3.1 from \cite{lmt1n}.

\begin{ex}\label{ex0}
\normalfont
Let  $\mu$ be an arbitrary probability measure on $\mathfrak B(G)$ and let $Q_{n}(\theta)$ be probability measures on $\mathfrak B(\varUpsilon)$ 
for all $n\in\mathbb N$ and for any fixed $\theta\in G$, which is absolutely continuous with respect to Lebesgue measure $\lambda$ on $\mathfrak B$. Moreover, suppose that there exists a $\B(G)$-measurable function $h:G\longmapsto \R$ such that $Q_{n}(\theta)={\bf{K}}\left(h(\theta)\right)$ for any $n\in\N$ and $\theta\in G$, where for any $B\in \mathfrak B(\varUpsilon)$ the function ${\bf{K}}\left(h(\bullet)\right)(B):G\longmapsto\R$ is $\B(G)$-measurable. Put $\wt{P}_{\theta}:=\otimes_{n\in\mathbb N} Q_{n}(\theta)$ for any $\theta\in G$.
Define the set-function $P(E) :=\int  \wt{P}_{\theta}(E^{\theta})\mu(d\theta)$, for each $E\in\vS$, where $E^{\theta}$ is the $\theta$-section of $E$, and put $P_{\theta}:=\wt{P}_{\theta}\otimes \delta_{\theta}$ for any $\theta\in G$, where $\delta_\theta$ is the Dirac measure at $\theta$. Then $P$ is a probability measure on $\vS$ and $\{P_{\theta}\}_{\theta\in G}$ is a disintegration of $P$ over $\mu$ consistent with the canonical projection $\pi_{G}$ from $\vO$ onto $G$ (compare  \cite{lmt1n}, Example 3.1).

Clearly, putting $\vT:= \pi_{G}$ we get $P_{\vT}=\mu$. Set $W_{n}:= \pi_{n}$ for any $n \in\mathbb N$, where $\pi_{n} :\vO\longmapsto\vY$ is the canonical projection, and $W:=\{W_{n}\}_{n\in\N}$. Put $T_n:=\sum_{k=1}^{n} W_k$ for any $n\in\N_0$ and $T:=\{T_n\}_{n\in\N_0}$, and let $N:=\{N_t\}_{t\in\R_+}$ be the counting process induced by $T$ by means of $
N_t:=\sum_{n=1}^{\infty}\chi_{\{T_n\leq t\}}$ for all $t\in\R_+$ (cf. e.g \cite{Sc}, Theorem 2.1.1). Applying the same arguments as in \cite{lmt1n}, Example 3.1, we get that  $N$  is a $P$-eMRP$({\bf K}(h(\vT)))$.

\end{ex}

In the next example the real-valued random variable $\wh{\vT}$ is distributed according to the Gamma law, a common choice in Risk Theory.

\begin{ex}\label{ex1}
\normalfont
Let $G:=\vY$, let $\xi={\bf{IGa}}(\alpha,\beta)$, with $\alpha,\beta>0$ be a probability measure on $\mathfrak B(\vY)$ i.e.
$$
\xi(B):=\int_B \frac{\beta^{\alpha}}{\Gamma(\alpha)}\cdot t^{-\alpha-1}\cdot e^{-\frac{\beta}{t}}\cdot\chi_{(0,\infty)}(t)\,\lambda(dt)
\quad\mbox{for each}\;\;B\in\mathfrak B(\vY)
$$
and let $h:\vY\longmapsto \R$ be defined by  $h(\theta):=\frac{1}{\theta}$ for any $\theta\in\vY$. Fix now on arbitrary $\theta\in\vY$ and define the probability measures $Q_{n}(\theta)$ by means of $Q_{n}(\theta):={\bf{Exp}}(h(\theta))$ for all $n\in\N$. Let $(\vO,\vS,P)$, $\vT$, $N$, $W$ and $\{P_\theta\}_{\theta\in\vY}$ be as in Example \ref{ex0} with $G=\vY$ and $\xi$ in the place of  $\mu$. 

Define the map $C\in\mathcal L^1(P_{h(\vT)})$ by $C(h(\theta)):=h(\theta)$ for any $\theta\in\vY$, and for any fixed $\theta\in\vY$ define the density $f_{h(\theta)}:=F'_{h(\theta)}$  by $f_{h(\theta)}(t):=h(\theta)\cdot e^{-h(\theta)t}$ for any $t>0$. Clearly, for any fixed $\theta\in\vY$, the density $f_{h(\theta)}$ is dominated by $C(h(\theta))$, and the function $p_h:\vY\longmapsto \vY$ defined by means of $p_h(\theta):=\lim_{t\rightarrow 0}f_{h(\theta)}(t)=h(\theta)$ for any $\theta\in\vY$, is positive and injective; hence $\{P_\theta\}_{\theta\in\vY}$, $N$ and $h$ satisfy Assumption \ref{as}. 

Let $\wh{\vT}:=h\circ\vT$ and put $Q_{\wh{\theta}}(E):=\left(P_{\bullet}(E)\circ h^{-1}\right)(\wh{\theta})$ for any $\wh{\theta}>0$ and $E\in\vS$. Then $\{Q_{\wh{\theta}}\}_{\wh{\theta}>0}$ is a disintegration of $P$ over $P_{\wh{\vT}}$ consistent with $\wh{\vT}$, condition $\left(Q_{\wh{\theta}}\right)_{W_n}={\bf{Exp}}(\wh{\theta})$ holds true for any $n\in\N$ and $\wh{\theta}>0$, and the process $W$ is $Q_{\wh{\theta}}$-independent  (see \cite{lmt1n} Lemma 2.4). Thus due to \cite{lm1v}, Proposition 4.4, we obtain  that $N$ is a $P$-MPP$(\wh{\vT})$.

Clearly, all  assumptions of Theorems \ref{thm1} and \ref{thm2} are satisfied and so are their conclusions. In particular, each of its assertions $(i)$ to $(iv)$ is valid.
\end{ex}

In our next  example the real-valued random variable $\wh{\vT}$ is distributed according to the Lognormal law, a common choice in Reliability Theory.
\medskip

\begin{ex}\label{ex2}
\normalfont
Let $G:=\R$, let $\rho=\mathbf{N}(\mu,\sigma^2)$, with $(\mu,\sigma^2)\in \R\times (0,\infty)$ be a probability measure on $\mathfrak B$ i.e.
$$
\rho(B):=\int_B \frac{1}{\sigma\sqrt{2\pi}}\,\cdot e^{-\frac{(t - \mu)^2}{2 \sigma^2}}\cdot\chi_{\R}(t)\;\lambda(dt)
\quad\mbox{for any}\;\;B\in\B
$$
and let $h:\R\longmapsto\R$ be defined by  $h(\theta):=e^{\theta}$ for any $\theta\in\R$.  Fix  on arbitrary $\theta\in\R$ and define the probability measures $Q_{n}(\theta)$ by means of $Q_{n}(\theta)\abc{:}={\bf{Exp}}(h(\theta))$ for all $n\in\N$. Let $(\vO,\vS,P)$, $\vT$, $N$, $W$ and $\{P_\theta\}_{\theta\in\R}$ be as in Example \ref{ex0} with $G=\R$ and $\rho$ in the place of $\mu$. 

Define the map $C\in\mathcal L^1(P_{h(\vT)})$ by $C(h(\theta)):=h(\theta)$ for any $\theta\in\R$, and for any fixed $\theta\in\R$ define the density $f_{h(\theta)}:=F'_{h(\theta)}$  by $f_{h(\theta)}(t):=h(\theta)\cdot e^{-h(\theta)t}$ for any $t>0$. Clearly, for any fixed $\theta\in\R$, the density $f_{h(\theta)}$ is dominated by $C(h(\theta))$, and the function $p_h:\R\longmapsto \vY$ defined by means of $p_h(\theta):=\lim_{t\rightarrow 0}f_{h(\theta)}(t)=h(\theta)$ for any $\theta\in\R$, is positive and injective; hence $\{P_\theta\}_{\theta\in\R}$, $N$ and $h$ satisfy Assumption \ref{as}.

Let $\wh{\vT}:=h\circ\vT$ and put $Q_{\wh{\theta}}(E):=\left(P_{\bullet}(E)\circ h^{-1}\right)(\wh{\theta})$ for any $\wh{\theta}>0$ and $E\in\vS$. 
It then follows as in Example \ref{ex1} that there exists a MPP$(\wh{\vT})$ with a lognormally distributed real-valued random variable $\wh{\vT}$ and that all  assumptions of Theorem \ref{thm2} are satisfied and so are its conclusions. In particular, each of its assertions $(i)$ to $(iv)$ is valid.
\end{ex}


\section{Counter-examples}\label{cexs}

The next counter-examples show that there exist probability spaces and counting processes on them satisfying  assertions \textit{(i)}, \textit{(ii)} and \textit{(iv)} but not assertion \textit{(iii)}  of Theorems \ref{thm1} and \ref{thm2}.

Moreover, the assumptions of  Theorem \ref{thm1}, concerning the perfectness of the measure $P$ and the countability of $\vS$, are not  valid, showing in this way that they are essential for the equivalence \textit{(i)} $\Longleftrightarrow$ \textit{(iii)}. The same examples show that  the assumption of Theorem \ref{thm2} concerning the existence of a disintegration consistent with $\vT$ is not valid; hence it is essential for the equivalence of \textit{(iii)} with any of the assertions \textit{(i)}, \textit{(ii)}, \textit{(iv)}.\smallskip

To present our counter-examples we need  the following result 

\begin{lem}\label{lem6}
Let $B$ be a subset of $\vO$ with $P^\ast(B)=1$ and $P_\ast(B)=0$. Put $\vS_b:=\sigma(\vS\cup\{B\})$ and define $R:=P_b:\vS_b\longmapsto[0,1]$ by means of $R(\wt{D}):=P^{\ast}(\wt{D}\cap B)$ for any 
$\wt{D}\in\vS_b$. Then there does not exist any d-dimensional random vector $\Psi$ on $\vO$  such that there exists a disintegration $\{R_{\psi}\}_{\psi\in\R^d}$ of $R$ over $R_{\Psi}$ consistent with $\Psi$. In particular, if $\vS$ is countably generated then $R$ is non-perfect.
\end{lem}

{\bf Proof.} Assume, if possible, that there exists a $d$-dimensional random vector $\Psi$ on $\vO$ such that there exists a disintegration 
$\{R_\psi\}_{\psi\in\R^d}$ of $R$ over $R_\Psi$ consistent with $\Psi$. For any $\omega\in\vO$ put  
$$
Q_\omega(E):=R_{\Psi(\omega)}(E)
$$
for each $E\in\wt{\vS}$.
\smallskip

\begin{it}
Claim 1. The family $\{Q_\omega\}_{\omega\in\vO}$ is a subfield regular conditional probability for $R$ over 
$R\upharpoonright \mathcal F$ with $\F=\sigma(\Psi)$.
\medskip

Proof.
\end{it} 
For the definition of a subfield regular conditional probability see \cite{fa}, Section 2. Clearly for any fixed $\omega\in\vO$ the set-function $Q_\omega$ is a probability measure on $\vS_b$, and for any fixed  $F\in\F$ the function $\omega\longmapsto Q_\omega(F)$ is $\F$-measurable. Furthermore, for each $F\in\F$ and $E\in\vS_b$ we have

\begin{eqnarray*}
\int_F Q_\omega(E)\, R(d\omega)&=& 
\int_{F} R_{\Psi(\omega)} (E)\,R(d\omega)
= \int_{F} \E_{R}[\chi_E\mid\F](\omega)\,R(d\omega)\\
&=&\int_F \chi_E\, dR=\int \chi_{F\cap E}\, dR;
\end{eqnarray*}
hence 
\begin{equation}
\label{1a}
\int_F Q_\omega(E)\, R(d\omega)=R(E\cap F)
\end{equation}
where the second equality follows from the assumption that the restriction of $\{R_\psi\}_{\psi\in\R^d}$ is a disintegration of $R$ over $R_\Psi$ consistent with $\Psi$, together with \cite{lm1v}, Lemma 2.5 (i). As a consequence, we get that 
$\{Q_\omega\}_{\omega\in\vO}$ is a subfield regular conditional probability for $\wt P$ over $\wt P\uph\F$. This completes the proof of Claim 1.\hfill$\Box$ \medskip

\begin{it}
Claim 2. There exists a $P$-null set $N\in\mathcal F$ such that for each $A\in\mathcal F$ condition $Q_\omega(A)=1$ holds true for any 
$\omega\in N^c\cap A$.
\medskip

Proof.
\end{it} 
Since $\{R_\psi\}_{\psi\in\R^d}$ is a  disintegration of $R$ over $R_\Psi$  consistent with $\Psi$ we get

\begin{equation}
\forall\; A\in\F\quad \exists\; N_A\in\F_0\;\; \forall\; \omega\in A\cap N^c_A\qquad [Q_\omega(A)=1],
\label{eq:1}
\end{equation}

where $\F_0$ is the set of all $P$-null sets in $\F$. Notice that $\F$ is countably generated since $\mathfrak B_d$ is so. Let $\mathcal G$ be a countable generator of $\F$. Without loss of generality we may and do assume that $\mathcal G$ is closed under finite intersections.
Since $\mathcal G$ is a countable generator of $\F$ condition (\ref{eq:1}) can be rewritten as 
$$
\forall\;n\in\N\;\;\forall\;A_n\in\mathcal G\quad \exists\; N_{A_n}\in\F_0\;\;\forall\;\omega\in A_n\cap N^c_{A_n}\qquad[Q_\omega(A_n)=1].
$$

So, setting $N:=\bigcup_{n\in\N} A_n$ we get that $N\in\F$ with $P(N)=0$.
Let us denote now by $\mathcal D$ the class of all sets $A\in\F$ such that $Q_{\omega}(A)=1$ for each $\omega\in N^c\cap A$. Then 
it can be easily proven that $\mathcal D$ is a Dynkin class; hence  by 
the Monotone Class Theorem the claim  follows.
\hfill$\Box$\medskip

But by condition \eqref{1a} and our assumption that $P^\ast(B)=1$,  we get for every $F\in\mathcal F$ that
$$
\int_F Q_\omega(B)\;R(d\omega)=R(F\cap B)=R(F)=\int_F \chi_F(\omega)\;R(d\omega),
$$
implying that $R(D)=0$, where $D:=\{\omega\in\vO:Q_\omega(B)\neq 1\}$

Put $E:=D\cup N$. For any $\omega\in E^c$ we get $Q_\omega(\{\omega\})=1$ and $Q_\omega(B)=1$; hence $Q_\omega(B\cap\{\omega\})=1$, implying $B\cap\{\omega\}\neq\emptyset$ or $\omega\in B$. Thus we get $E^c\subseteq B$ or equivalently $B^c\subseteq E$, 
implying $1=P^\ast(B^c)\leq R(E)$; hence $R(E)=1$, a contradiction. 

In particular, if $\vS$ is countably generated then $\vS_b$ is so; hence applying \cite{fa}, Theorem 6, we deduce that $R$ is non-perfect.
\hfill$\Box$

\begin{rem}\label{Bern}
\normalfont
Let $\vO$ be an uncountable Polish space and $P$ a non-atomic Borel measure on $\vS:=\B(\vO)$. It should be  known that there always exists a set $B\subseteq\vO$ such that $P^\ast(B)=1$ and $P_\ast(B)=0$. But since we could not find it in the literature, we insert a short proof for completeness sake: Let $(\vO, \wh{\vS}, \wh{P})$ be the completion of $(\vO, \vS, P)$. Then $(\vO, \wh{\vS}, \wh{P})$ is isomorphic to the 
Lebesgue probability space $([0,1], \mathcal L([0,1]), \lambda)$ (cf. e.g. \cite{fr3}, Corollary 344K). By \cite{co}, Proposition 1.4.11, there exists a subset $\wt{A}
$ of $\R$ such that each Lebasgue measurable set that is included in $\wt{A}$ or $\wt{A}^c$ is a $\lambda$-null set. Put $A:=[0,1]\cap \wt{A}$ and $A_1:=[0,1]\cap \wt{A}^c$. 

The subsets $A$ and $A_1$ cannot be both Lebesgue measurable, since if they were so, then we would get that $\lambda(A)=\lambda(A_1)=0$, implying  $0=\lambda(A\cup A_1)=\lambda([0,1])=1$, a contradiction. Thus, if $A$ is non-Lebesgue measurable, we infer that $\lambda_\ast(A)=0$ and $\lambda^\ast(A)=1$. Without loss of generality we may and do assume that $A$ is non-Lebesgue measurable. So, letting $f:[0,1]\longmapsto \vO$ be an isomorphism between the Lebesgue probability space on $[0,1]$ and $(\vO,\wh{\vS},\wh{P})$, we get that $B:=f(A)$ is the desired set.
\end{rem}

\begin{rem}\label{null}
\normalfont
Let $P$, $\vT$ and $\{P_\theta\}_{\theta\in\vY}$ be as in Example \ref{ex1}. Fix on arbitrary $\theta\in\vY$ and put $\vS_0:=\{L\in\vS : P(L)=0\}$ and $\vS_{0,\theta}:=\{L\in\vS : P_\theta(L)=0\}$. \smallskip

{\bf(a)} Since for any fixed $E\in\wt\vS$ the function $\theta\longmapsto \wt P_\theta(E):=\otimes_{n\in\N}{\bf Exp}(h(\theta))(E)$ is continuous, it can be easily seen that $\vS_0=\vS_{0,\theta}$ for every $\theta\in\vY$, implying that $P^\ast(B)=P^\ast_\theta(B)=1$ and $P_\ast(B)=(P_\theta)_\ast(B)=0$. Thus, the probability measure $P_\theta$ can be extended to the probability measure $P_{\theta,b}:\vS_b\longmapsto [0,1]$, defined by means of $P_{\theta,b}(\wt D):=P^\ast_\theta(\wt D\cap B)$ for any $\wt D\in\vS_b$. Then for any fixed $D\in\vS_b$ the function $\theta\longmapsto P_{\theta,b}(\wt D)$ is $\B(\vY)$-measurable.

{\bf(b)} For any $\theta\in\vY$ denote by $\wh \vS_\theta$ the completion of $\vS$ with respect to $P_\theta$. It then follows that $\wh\vS=\wh\vS_\theta$ for any $\theta\in\vY$; hence  each completed probability measure $\wh P_\theta$ is defined on $\wh\vS$ and for any $E\in\wh\vS$ the function $\theta\longmapsto \wh P_{\theta}(E)$ is $\B(\vY)$-measurable.
\end{rem}

\begin{ex}\label{ex6}
\normalfont
Let $(\vO,\vS,P)$, $N$, $\vT$, $h$ and $\wh\vT$ be as in Example \ref{ex1}. Then all the assumptions of Theorems \ref{thm1} and \ref{thm2} 
are satisfied, and so the equivalence of all assertions (under $P$) follows for each of the above two theorems. In particular, recall that each assertion of both theorems is valid.

Since by  construction $(\vO,\vS,P)$ is an uncountable non-atomic Polish probability space, it follows by Remark \ref{Bern} that there exists a set $B\subseteq \vO$ such that  $P^\ast(B)=1$ and $P_\ast(B)=0$; hence we may define  $R$ and $\vS_b$ as in Lemma \ref{lem6}. It can be easily seen that assertions $(i)$ and $(iv)$ of Theorems \ref{thm1} and \ref{thm2} remain valid under $R$, and taking into account Remark \ref{null} so does assertion $(ii)$. In particular, assertions $(i)$, $(ii)$ and $(iv)$ of Theorem \ref{thm2} are equivalent under $R$.

According to Lemma \ref{lem6}, there does not exist any real-valued random variable $\Psi$ on $\vO$ such that there exists a disintegration $\{R_\psi\}_{\psi\in\R}$ of $R$ over $R_\Psi$ consistent with $\Psi$, implying that assertion $(iii)$ of Theorems \ref{thm1} and \ref{thm2} fails.

Note that, due to Lemma \ref{lem6}, neither the part of the Assumption \ref{as1} concerning the existence of a real-valued random variable $\vT$ on $\vO$ and a disintegration of $R:=P_b$ over $R_\vT$ consistent with $\vT$ nor the perfectness assumption of Theorem \ref{thm1} for the probability space $(\vO,\vS_b,R)$ hold true. Thus, both assumptions are not necessary for the equivalences $(i)\Longleftrightarrow (ii) \Longleftrightarrow (iv)$ and the perfectness assumption is essential for the equivalence $(i)\Longleftrightarrow (iii)$.\smallskip

Concerning Theorem \ref{thm2}, due again to Lemma \ref{lem6}, the assumption of the existence of a disintegration of $R$ over $R_\vT$ consistent with $\vT$ is not necessary for the equivalences $(i)\Longleftrightarrow (ii) \Longleftrightarrow (iv)$ but it is essential for the equivalence $(i)\Longleftrightarrow (iii)$.
\end{ex}

\begin{ex}\label{ex7}
\normalfont
Let $(\vO,\vS,P)$, $N$, $\vT$, $h$ and $\wh\vT$ be as in Example \ref{ex1}. Then all the assumptions of Theorems \ref{thm1} and \ref{thm2} are satisfied, and so the equivalence of all assertions (under $P$) follows for each of the above two theorems. In particular, recall that each assertion of both theorems is valid. Let $(\vO,\wh{\vS},\wh{P})$ be the completion of $(\vO,\vS,P)$. It can be easily seen that assertions $(i)$ and $(iv)$ of Theorems \ref{thm1} and \ref{thm2} remain valid under $\wh{P}$, and taking into account Remark \ref{null} so does assertion $(ii)$. In particular, assertions $(i)$, $(ii)$ and $(iv)$ of Theorem \ref{thm2} are equivalent under $\wh{P}$.

First notice that the probability measure $\wh{P}$ is perfect since $P$ is (cf. e.g. \cite{fr4}, Proposition 451G(c)(i)), but  $\wh{\vS}$ is not countably generated; hence the countability assumption in Theorem \ref{thm1} fails.\medskip

\begin{it}
Claim. There does not exist any real-valued random variable $\Psi$ on $\vO$ such that there exists a disintegration $\{Z_{\psi}\}_{\psi\in\R}$ of $\wh{P}$ over $\wh{P}_{\Psi}$ consistent with $\Psi$.\medskip

Proof.
\end{it}
Assume, if possible, that there exists  $\Psi$ on $\vO$ such that there exists a disintegration $\{Z_{\psi}\}_{\psi\in\R}$ of $\wh{P}$ over $\wh{P}_{\Psi}$ consistent with $\Psi$. Fix on an arbitrary $A\in\wh{\vS}$ and define the function  $S_\bullet(A):\vO\longmapsto[0,1]$ by 
means of 
$$
(S_\bullet(A))(\omega):=S_\omega(A):=(Z_\bullet(A)\circ\Psi)(\omega)\;\;\text{for every}\;\; \omega\in\vO.
$$ 
Then by the same arguments as in the proof of Lemma \ref{lem6} we get that $\{S_\omega\}_{\omega\in\vO}$ is a  subfield r.c.p. of $\wh{P}$ over $\wh{P}\upharpoonright\sigma(\Psi)$. Since $\sigma(\Psi)$ is countably generated, it follows as in Lemma \ref{lem6} that there exists a set $N\in\sigma(\Psi)$ such that $P(N)=0$ and for any $A\in\sigma(\Psi)$ condition $S_\omega(A)=1$ holds true for any $\omega\in N^c\cap A$.
Choose a set $D\subseteq N^c$ such that $D\notin\sigma(\Psi)$ but $D\in\wh{\vS}$. 
Such a choice is possible, since the cardinality of $\sigma(\Psi)$ is ${\bf c}$, ${\bf c}$ the cardinality of the continuum, while the cardinality of $\wh{\vS}$ is $2^{\bf c}$. Then for each $\omega\notin N$ we obtain
$$
1=S_\omega(\{\omega\})\leq S_\omega(D)\leq 1\qquad\text{if}\quad \omega\in D
$$
and
$$
1=S_\omega(\{\omega\})\leq S_\omega(D^c)\leq 1\qquad\text{if}\quad \omega\in D^c.
$$
Thus, $D=N^c\cap\{\omega\in\vO:S_\omega(D)=1\}\in\sigma(\Psi)$, a contradiction.
\hfill$\Box$ \medskip

 As a consequence, it follows that neither Assumption \ref{as1} of Theorem \ref{thm1} nor the assumption of Theorem \ref{thm2} concerning the existence of a disintegration hold true.
Moreover, the above claim  yields that assertion $(iii)$ of Theorems \ref{thm1} and \ref{thm2} fails; hence the countability assumption for $\vS$ is essential for the equivalence of $(i)$ and $(iii)$. 
\end{ex}

\begin{rem}\label{rem3}
\normalfont
The above two counter-examples answer to the negative  \cite{lmt1n}, Question 2.14, concerning the necessity of the assumptions  of the existence of a disintegration of $P$ over $P_\vT$ consistent with $\vT$ of  Proposition 2.7 and Theorem 2.11 from \cite{lmt1n}, for the validity of their conclusions. \medskip

In fact, let $(\vO,\wh\vS,\wh P)$ and $N$, $h$, $\vT$, $\wh\vT$ and $\{P_\theta\}_{\theta>0}$ 
be as in Example \ref{ex7}, and let $\{\wh P_\theta\}_{\theta>0}$ be as in Remark \ref{null}. Then $N$ is $\wh P$-eMRP$({\bf K}(h(\vT)))$ satisfying together with $h$ and  $\{\wh P_\theta\}_{\theta>0}$ Assumption \ref{as}, but according to the claim  of Example \ref{ex7}, $\{\wh P_\theta\}_{\theta>0}$ can not be a disintegration of $\wh P$ over $\wh P_\vT$ consistent with $\vT$. For the counting process $N$ the following are equivalent
\begin{itemize}
\item[$(a)$] $N$ is a $\wh P$-MPP$(\wh\vT)$;
\item[$(b)$] $N$ has the $\wh P$-multinomial property;
\item[$(c)$] $N$ has the $\wh P$-Markov property.
\end{itemize}
For the definition of the multinomial property we refer to e.g. \cite{SZ}, Section 2 page 2.

Ad $(a)\Longleftrightarrow (b)$: $N$ is a $\wh P$-MPP$(\wh\vT)$ if and only if it is a  $\wh P$-MPP$(P_{\wh\vT})$ (see Example \ref{ex7}) if and only if it has the $\wh P$-multinomial property (see \cite{SZ}, Theorem 4.2).\smallskip 

Ad $(a)\Longleftrightarrow (c)$: $N$ is a $\wh P$-MPP$(\wh\vT)$   if and only if it is a $\wh P$-MPP$(\{\wh P_{\theta}\}_{\theta>0},\nu)$ (see Example \ref{ex7}). But the latter, due to  \cite{hu}, Theorem 3, is equivalent  with the $\wh P$-Markov property for $N$ (see also \cite{lmt1n}, Remark 2.8 (d) and Example 3.3).\smallskip

Thus, in the above set-up we constructed a probability space $(\vO,\wh\vS,\wh P)$, a family of probability measure $\{\wh P_\theta\}_{\theta>0}$ on $\wh\vS$, and a counting process $N$ being a $\wh P$-eMRP$({\bf K}(h(\vT)))$ satisfying together with $h$ and  $\{\wh P_\theta\}_{\theta>0}$ Assumption \ref{as} such that the conclusions of  Theorem 2.11 from \cite{lmt1n}, hold   true but $\{\wh P_\theta\}_{\theta>0}$ is not necessarily a disintegration of $\wh P$ over $\wh P_\vT$ consistent with $\vT$.
\end{rem}

\medskip

{\small
\noindent{\sc D.P. Lyberopoulos, N.D. Macheras and S.M Tzaninis}\\
{\sc Department of Statistics and Insurance Science}\\
{\sc University of Piraeus, 80 Karaoli and Dimitriou street}\\
{\sc 185 34 Piraeus, Greece}\\
E-mail: {\tt dilyber@webmail.unipi.gr},$\;$ {\tt macheras@unipi.gr}$\;$ {\sc and}$\;$ {\tt stzaninis@unipi.gr}
}
\end{document}